\documentclass[12pt]{amsart}

\usepackage{amssymb}
\usepackage{amsfonts}
\usepackage{amsmath}

\theoremstyle{plain}
\newtheorem{thm}{Theorem}[section]
\newtheorem{lem}{Lemma}[section]
\newtheorem{cor}{Corollary}[section]

\theoremstyle{definition}

\numberwithin{equation}{section}

\begin{document}
\baselineskip=17pt

\title{The Radius in Matrix Algebras---Examples and Remarks}

\author{Moshe Goldberg}

\address{Department of Mathematics,
Technion -- Israel Institute of Technology,
Haifa 32000, Israel}

\email{mg@technion.ac.il}



\begin{abstract}
The main purpose of this note is to illustrate how the {\em radius} in a finite-dimensional
power-associative algebra over a field $\mathbb{F}$, either $\mathbb{R}$ or $\mathbb{C}$, may
change when the multiplication in this algebra is modified. Our point of departure will be
$\mathbb{F}^{n \times n}$, the familiar algebra of $n \times n$ matrices over $\mathbb{F}$ with
the usual matrix operations, where it is known that the radius is the classical spectral radius.
We shall alter the multiplication in $\mathbb{F}^{n \times n}$ in three different ways and compute,
in each case, the radius in the resulting algebra.
\end{abstract}

\dedicatory{To my daughter Maya on her 40th birthday}

\maketitle

\section{\label{sec:1}Introduction: the radius and its basic properties}

Let $\mathcal{A}$ be a finite-dimensional algebra over a field $\mathbb{F}$, either $\mathbb{R}$ or
$\mathbb{C}$. We shall assume that $\mathcal{A}$ is {\em power-associative}, i.e., that the subalgebra
of $\mathcal{A}$ generated by any one element is associative; thus ensuring that {\em powers of each
element in $\mathcal{A}$ are unambiguously defined.}

As usual, by a {\em minimal polynomial} of an element $a$ in $\mathcal{A}$ we mean a monic polynomial
of lowest positive degree with coefficients in $\mathbb{F}$ that annihilates $a$.

With this familiar definition, we may cite:

\begin{thm}
[{[G1, Theorem 1.1]}]
\label{thm:1.1}
Let $\mathcal{A}$ be a finite-dimensional power-associative algebra over $\mathbb{F}$. Then:

{\em (a)} Every element $a\in \mathcal{A}$ possesses a unique minimal polynomial.

{\em (b)} The minimal polynomial of $a$ divides every other polynomial over $\mathbb{F}$ that
annihilates $a$.
\end{thm}

Denoting the minimal polynomial of an element $a \in \mathcal{A}$  by $p_a$, we follow [G1] and define
the {\em radius} of $a$ to be the nonnegative quantity
$$
r(a)=\max \{ |\lambda|:~\lambda \in \mathbb{C}, \lambda \textrm{ is a root of } p_a \}.
$$
The radius has been computed for elements in several well-known finite-dimensional power-associative
algebras. For instance, it was recently shown in [GL1] that the radius in the Cayley--Dickson algebras
is given by the corresponding Euclidean norm.

Another example emerged in [G1, page 4060], where it was established that if $\mathcal{A}$  is an
arbitrary finite-dimensional matrix algebra over $\mathbb{F}$  with the usual matrix operations, then
the radius of a matrix $A\in \mathcal{A}$ is given by the classical spectral radius,
$$
\rho(A)=\max \{ |\lambda|:~\lambda \in \mathbb{C}, \lambda \textrm{ is an eigenvalue of } A \}.
$$

With this last example in mind, we recall the following theorem which asserts that the radius retains
some of the most basic properties of the spectral radius not only in finite-dimensional matrix algebras
with the usual matrix operations, but in the
general finite-dimensional power-associative case as well.

\begin{thm}
[{[G1, Theorems 2.1 and 2.4]}]
\label{thm:1.2}
Let $\mathcal{A}$ be a finite-dimensional power-associative algebra over $\mathbb{F}$. Then:

{\em (a)} The radius $r$ is a nonnegative function on $\mathcal{A}$.

{\em (b)} The radius is homogeneous, i.e., for all $a\in \mathcal{A}$ and $\alpha \in \mathbb{F}$,
$$
r(\alpha a)=|\alpha| r(a).
$$

{\em (c)} For all $a \in \mathcal{A}$ and all positive integers $k$,
$$
r(a^k)=r(a)^k.
$$

{\em (d)} The radius vanishes only on nilpotent elements of $\mathcal{A}$.

{\em (e)} The radius is a continuous function on $\mathcal{A}$.\footnote{Naturally, a real-valued
function on a finite-dimensional algebra $\mathcal{A}$ is said to be {\em continuous} if it is
continuous with respect to the (unique) finite-dimensional topology on $\mathcal{A}$.}
\end{thm}

As a final introductory remark we mention that an analysis of the relevance of the radius to stability
of subnorms and to the Gelfand formula can be found in [G1], [G2], [G3], and [GL1].

\section{\label{sec:2}Examples of radii in matrix algebras}

Our main purpose in this note is to illustrate how the radius in a finite-dimensional power-associative
algebra may change when the multiplication in this algebra is modified. Selecting a positive integer
$n$, $n \geq 2$, our point of departure will be $\mathbb{F}^{n \times n}$, the familiar algebra of
$n \times n$ matrices over $\mathbb{F}$, either $\mathbb{R}$ or $\mathbb{C}$, with the usual matrix
operations. By what we already know about the radius in arbitrary finite-dimensional matrix algebras over $\mathbb{F}$ with the usual operations, we may register the following result which can also be derived directly from the fact that the roots of the minimal polynomial of a matrix $A$ in $\mathbb{F}^{n \times n}$ are the eigenvalues of $A$.

\begin{thm}
\label{thm:2.1}
The radius of a matrix $A$ in $\mathbb{F}^{n \times n}$ is given by
$$
r(A)= \rho(A),
$$
where $\rho$ denotes the spectral radius.
\end{thm}

The multiplication in $\mathbb{F}^{n \times n}$ can be altered, of course, in a myriad of ways.
Often, however, computing the radius in the newly obtained algebra will remain out of reach. In what
follows, we shall modify the multiplication in $\mathbb{F}^{n \times n}$ in three different ways, and
calculate the radius in each case.

We embark on our plan by replacing the standard multiplication in $\mathbb{F}^{n \times n}$ by the
well-known Hadamard product which, for any two $n \times n$  matrices $A=(\alpha_{ij})$  and
$B=(\beta_{ij})$, is defined entry-wise by
$$
A\circ B = (\alpha_{ij} \beta_{ij}).
$$

The resulting algebra, denoted by $\mathbb{F}_H^{n \times n}$, has been extensively studied in the
literature (see for example Chapter 5 in [HJ] and the references at the end of that chapter).
Obviously, $\mathbb{F}_H^{n \times n}$ is distributive, commutative, and associative; and its unit
element is given by $E$, the $n \times n$ matrix all of whose entries are 1.

Denoting the $k$-th power of a matrix $A=(\alpha_{ij})$ in $\mathbb{F}_H^{n \times n}$ by $A^{[k]}$,
we see that
\begin{equation}
\label{eq:2.1}
A^{[k]} = (a_{ij}^k), \quad k=1,2,3,\ldots.
\end{equation}
Assisted by this observation, we can now post:

\begin{thm}
\label{thm:2.2}
The radius of a matrix $A=(\alpha_{ij})$ in $\mathbb{F}_H^{n \times n}$ is given by the sup norm of $A$, i.e.,
$$
r(A)= \max_{i,j} |\alpha_{ij}|.
$$
\end{thm}

\noindent {\em Proof.}
Select a matrix $A = (\alpha_{ij})$ in $\mathbb{F}_H^{n \times n}$, and let $\zeta_1, \ldots, \zeta_s$
$(1 \leq s \leq n^2)$ be a list of all the distinct entries of $A$ (so that each $\alpha_{ij}$ equals precisely one of the $\zeta_l$'s). Let
$$
p_A(t) = t^m + \alpha_{m-1}t^{m-1} + \cdots + \alpha_1 t + \alpha_0
$$
be the minimal polynomial of $A$ in $\mathbb{F}_H^{n \times n}$, hence
$$
A^{[m]} + \alpha_{m-1}A^{[m-1]} + \cdots + \alpha_1 A^{[1]} + \alpha_0 E = 0.
$$
By (2.1), this can be equivalently written as
$$
\zeta_l^m + \alpha_{m-1}\zeta_l^{m-1} + \cdots + \alpha_1 \zeta_l + \alpha_0 = 0, \quad l = 1,\ldots,s.
$$
It follows that the $\zeta_l$  are roots of $p_A$; and since these roots are distinct, we infer that
the monic polynomial
$$
q(t) = (t - \zeta_1)(t - \zeta_2) \cdots (t - \zeta_s)
$$
must divide $p_A$. On the other hand, we notice that
$$
(A - \zeta_1 E) \circ (A - \zeta_2 E) \circ \cdots \circ (A - \zeta_s E) = 0;
$$
so $q$ annihilates $A$ in $\mathbb{F}_H^{n \times n}$. Appealing to Theorem 1.1(b), we conclude that $p_A$ must divide $q$; hence $p_A = q$, and the rest of the proof follows without difficulty.
\qed

Another way of altering the standard multiplication in $\mathbb{F}^{n \times n}$ is to replace it by
the familiar Jordan product
$$
A \cdot B = \frac{1}{2} (AB + BA),
$$
which turns $\mathbb{F}^{n \times n}$ into the {\em special Jordan algebra} $\mathbb{F}^{n \times n+}$
(e.g., [J, page 4, Definition 2]). Since $\mathbb{F}^{n \times n}$ is distributive, so is
$\mathbb{F}^{n \times n+}$. Further, both $\mathbb{F}^{n \times n+}$ and $\mathbb{F}^{n \times n}$ share
the same unit element, the $n \times n$ identity matrix $I$. We observe, however,
that $\mathbb{F}^{n \times n+}$, unlike $\mathbb{F}^{n \times n}$, is commutative. Moreover, in contrast
with $\mathbb{F}^{n \times n}$, the algebra $\mathbb{F}^{n \times n+}$ is not associative, nor even
{\em alternative}.\footnote{As usual, we call an algebra $\mathcal{A}$ {\em alternative} if the subalgebra
generated by any two elements of $\mathcal{A}$ is associative.} Indeed, consider the matrices
$$
A = \begin{pmatrix}
  0 & 1 \\
  0 & 0 \\
\end{pmatrix} \oplus O_{n-2},
\quad
B = \begin{pmatrix}
  0 & 0 \\
  1 & 0 \\
\end{pmatrix} \oplus O_{n-2},
$$

where $O_{n-2}$ is the $(n-2) \times (n-2)$ zero matrix. Then,
$$
(A \cdot B) \cdot B = \frac{1}{4} (AB^2+2BAB+B^2A) = \frac{1}{2} BAB
= \frac{1}{2} B \neq 0 = A \cdot (B \cdot B),
$$
and alternativity is shattered.

Despite the fact that $\mathbb{F}^{n \times n+}$ is not alternative, it is power-associative. This is
so because powers of matrices in $\mathbb{F}^{n \times n+}$ coincide with those in
$\mathbb{F}^{n \times n}$, and hence are uniquely defined.

Turning to compute the radius in $\mathbb{F}^{n \times n+}$, we realize that it is a simple task:
Since $\mathbb{F}^{n \times n+}$ and $\mathbb{F}^{n \times n}$ have an identical linear structure,
and since raising to powers in $\mathbb{F}^{n \times n+}$ and $\mathbb{F}^{n \times n}$ coincide,
the minimal polynomials of a matrix $A$ in $\mathbb{F}^{n \times n+}$ and in $\mathbb{F}^{n \times n}$
are one and the same; thus, the radii of $A$ in $\mathbb{F}^{n \times n+}$ and in
$\mathbb{F}^{n \times n}$ come to the same thing, yielding:

\begin{thm}
\label{thm:2.3}
The radius of a matrix $A$ in $\mathbb{F}^{n \times n+}$ is given by
$$
r(A)= \rho(A).
$$
\end{thm}

This result is of particular interest, precisely because it tells us that altering the multiplication in
$\mathbb{F}^{n \times n}$ does not necessarily result in a different radius.

In our last example, we shall modify the multiplication in $\mathbb{F}^{n \times n}$ in a more
intricate way, by introducing the product
$$
A \ast B = (A'B')',
$$
where $A'$ is the matrix obtained from $A$ by replacing $\alpha_{1n}$, the $(1,n)$ entry of $A$, by its
negative, and where $A'B'$ is the usual product of $A'$ and $B'$ in $\mathbb{F}^{n \times n}$.

Denoting our new algebra by $\mathbb{F}^{n \times n}_\ast$, we remark that it is distributive
and associative. Indeed, for all $A,B,C \in \mathbb{F}^{n \times n}_\ast$  we have,
$$
A\ast(B+C) = (A'(B+C)')' = (A'B'+A'C')' = (A'B')'+(A'C')' =  A\ast B + A\ast C
$$
and similarly,
$$
(A+B)\ast C = A\ast C + B\ast C;
$$
so the distributive laws are in the bag. Furthermore, since
$$
(A\ast B)\ast C = (A'B')'\ast C = ((A'B')C')' = (A'(B'C'))' = A\ast (B'C')' = A\ast (B\ast C),
$$
associativity holds as well.

We also observe that the identity matrix $I$ constitutes the unit element in
$\mathbb{F}^{n \times n}_\ast$, since for all $A \in \mathbb{F}^{n \times n}_\ast$,
$$
A\ast I = (A'I')' = (A'I)' = A,
$$
and analogously, $I\ast A = A$. Lastly, we note that since $\mathbb{F}^{n \times n}$ is not
commutative, neither is $\mathbb{F}^{n \times n}_\ast$.

It seems interesting to mention that the algebra
$\mathbb{F}^{n \times n}_\ast$ possesses certain exotic properties which are not shared by either
$\mathbb{F}^{n \times n}$, $\mathbb{F}^{n \times n}_H$, or $\mathbb{F}^{n \times n+}$. For instance,
$\mathbb{F}^{n \times n}_\ast$ contains nilpotent matrices which have nonzero eigenvalues. To
substantiate this statement, let $A = (\alpha_{ij})$ be the $n\times n$ matrix all of whose entries
are zero except for $\alpha_{11}$, $\alpha_{1n}$, $\alpha_{n1}$, and $\alpha_{nn}$ which are given by
1, $-i$, $i$, and -1, respectively. It is not hard to verify that $A\ast A$ = 0, so $A$ is a nilpotent
matrix of index 2 in $\mathbb{F}^{n \times n}_\ast$. At the same time, we have
$$
\det(tI - A) = t^{n-2}(t^2 -2),
$$
so $\sqrt{2}$ and $-\sqrt{2}$ are eigenvalues of $A$.

Another property of $\mathbb{F}^{n \times n}_\ast$ which is not shared by our previous matrix algebras
lies in the fact that $\mathbb{F}^{n \times n}_\ast$ admits positive matrices whose squares are
negative.\footnote {By a {\em positive} matrix we mean here a nonzero matrix all of whose entries are
nonnegative. Similarly, a {\em negative} matrix is a nonzero matrix whose entries are all non-positive.}
For example, consider the $n \times n$  matrix $A = (\alpha_{ij})$ where $\alpha_{1n} = \alpha_{n1} =1$
and the rest of the entries vanish. While $A$ is positive, its squaring in
$\mathbb{F}^{n \times n}_\ast$ provides the negative matrix all of whose entries are zero, except for the
first and last entries along its diagonal which equal -1.

Turning to compute the radius in $\mathbb{F}^{n \times n}_\ast$, we denote the $k$-th power of a matrix
$A$ in this algebra by $A^{\langle k \rangle }$, and offer the following elementary observation.

\begin{lem}
\label{lem:2.1}
If $A \in \mathbb{F}^{n \times n}_\ast$, then
\begin{equation}
\label{eq:2.2}
A^{\langle k \rangle} = ((A')^k)', \quad k = 1,2,3, \ldots,
\end{equation}
where $(A')^{{k}}$ is the usual k-th power of $A'$ in $\mathbb{F}^{n \times n}$.
\end{lem}

\noindent {\em Proof.}
For $k = 1$ the assertion is trivial. So assuming (2.2) for $k$, we get
$$
A^{\langle k+1 \rangle} = A\ast A^{\langle k \rangle} = A\ast ((A')^k)' = (A'(A')^k)' = ((A')^{k+1})',
$$
and we are done.
\qed

With the above lemma in our grip, we may now proceed to record:

\begin{thm}
\label{thm:2.4}
The minimal polynomial of a matrix A in $\mathbb{F}^{n \times n}_\ast$ coincides with the minimal
polynomial of $A'$ in $\mathbb{F}^{n \times n}$.
\end{thm}

\noindent {\em Proof.}
Let
$$
p(t) = \alpha_m t^m + \cdots + \alpha_1 t + \alpha_0
$$
be a polynomial over $\mathbb{F}$ that annihilates $A$ in $\mathbb{F}^{n \times n}_\ast$; that is,
$$
\alpha_m A^{\langle m \rangle} + \cdots + \alpha_1 A^{\langle 1 \rangle} + \alpha_0 I = 0.
$$
By (2.2), this is equivalent to
$$
\alpha_m ((A')^m)' + \cdots + \alpha_1 (A')' + \alpha_0 I' = 0;
$$
or in other words, to
$$
\alpha_m (A')^m + \cdots + \alpha_1 A' + \alpha_0 I = 0.
$$
It follows that $p$ annihilates $A$ in $\mathbb{F}^{n \times n}_\ast$ if and only if $p$
annihilates $A'$ in $\mathbb{F}^{n \times n}$; so aided by Theorem 1.1, the proof follows.
\qed

An immediate consequence of Theorem 2.4 reads:

\begin{cor}
\label{cor:2.1}
The radii of $A$ in $\mathbb{F}^{n \times n}_\ast$ and of $A'$ in $\mathbb{F}^{n \times n}$ coincide.
\end{cor}

Finally, since the radius in $\mathbb{F}^{n \times n}$ is the spectral radius, we get:

\begin{thm}
\label{thm:2.5}
The radius of a matrix $A$ in $\mathbb{F}^{n \times n}_\ast$ is given by
$$
r(A) = \rho (A').
$$
\end{thm}

We conclude this note by pointing out that all our findings regarding $\mathbb{F}^{n \times n}_\ast$
hold verbatim when the product is defined by
$$
A\ast B = (A'B')'
$$
where now, $A'$ is obtained from $A$ by negating $\alpha_{n1}$, the $(n,1)$ entry of $A$.

The author is truly grateful to Thomas Laffey for helpful discussions.


\begin{thebibliography}{-------}

\bibitem[G1]{Goldberg 1}Moshe~Goldberg,
{\em Minimal polynomials and radii of elements in finite-dimensional power-associative algebras},
Trans. Amer. Math. Soc. $\mathbf{359}$ (2007), no. 8, 4055--4072.

\bibitem[G2]{Goldberg 2}Moshe~Goldberg,
\emph{Radii and subnorms on finite-dimensional power-associative algebras},
Linear Multilinear Algebra, $\mathbf{55}$ (2007), no. 5, 405--415.

\bibitem[G3]{Goldberg 3}Moshe~Goldberg,
{\em Stable subnorms on finite-dimensional power-associative algebras},
Electron. J. Linear Algebra $\mathbf{17}$ (2008), 359--375.

\bibitem[GL1]{Goldberg-Laffey}Moshe~Goldberg and Thomas J. Laffey,
{\em On the radius in Cayley--Dickson algebras},
Proc. Amer. Math. Soc. $\mathbf{143}$, (2015), no. 11, 4733--4744.

\bibitem[HJ]{Horn-Johnson}Roger~A.~Horn and Charles~R.~Johnson,
{\em Topics in Matrix Analysis}, Cambridge Univ. Press, Cambridge 1991.

\bibitem[J]{Jacobson}Nathan Jacobson,
{\em Structure and representations of Jordan algebras}, Amer. Math. Soc. Colloquium Publications,
Vol. XXXIX, Amer. Math. Soc., Providence, R.I. 1968.

\end{thebibliography}
\end{document}